# MINIMAX ESTIMATION WITH THRESHOLDING AND ITS APPLICATION TO WAVELET ANALYSIS

By Harrison H. Zhou[1] and J. T. Gene Hwang[2]

*Yale University and Cornell University*

Many statistical practices involve choosing between a full model and reduced models where some coefficients are reduced to zero. Data were used to select a model with estimated coefficients. Is it possible to do so and still come up with an estimator always better than the traditional estimator based on the full model? The James–Stein estimator is such an estimator, having a property called minimaxity. However, the estimator considers only one reduced model, namely the origin. Hence it reduces no coefficient estimator to zero or every coefficient estimator to zero. In many applications including wavelet analysis, what should be more desirable is to reduce to zero only the estimators smaller than a threshold, called thresholding in this paper. Is it possible to construct this kind of estimators which are minimax?

In this paper, we construct such minimax estimators which perform thresholding. We apply our recommended estimator to the wavelet analysis and show that it performs the best among the well-known estimators aiming simultaneously at estimation and model selection. Some of our estimators are also shown to be asymptotically optimal.

**1. Introduction.** In virtually all statistical activities, one constructs a model to summarize the data. Not only could the model provide a good and effective way of summarizing the data, the model if correct often provides more accurate prediction. This point has been argued forcefully in Gauch (1993). Is there a way to use the data to select a reduced model so that if the reduced model is correct, the model-based estimator will improve on the naive estimator (constructed using a full model) and yet never do worse than the naive estimator even if the full model is actually the only correct model? James–Stein estimation (1961) provides such a striking result under the normality assumption. Any estimator such as the James–Stein estimator that does no worse than the naive estimator is said to be minimax. See the precise discussion right before Lemma 1 of Section 2. The problem with the James–Stein positive part estimator is, however, that it selects only between two models: the origin and the full model. It is possible to construct estimators similar to the James–Stein positive part to select between the full model and another linear subspace. However, it always chooses between the two. The nice idea of George (1986a, b) in multiple shrinkage does allow the data to choose among several models; it, however, does not do thresholding, as is the aim of the paper.

Models based on wavelets are very important in many statistical applications. Using these models involves model selection among the full model or the models with smaller dimensions where some of the wavelet coefficients are zero. Is there a way to select a reduced model so that the estimator based on it does no worse in any case than the naive estimator based on the full model, but improves substantially upon the naive estimator when the reduced model is correct? Again, the James–Stein estimator provides such a solution. However, it selects either the origin or the full model. Furthermore, the ideal estimator should do thresholding; namely, it gives zero as an estimate for the components which are smaller than a threshold, and preserves (or shrinks) the other components. However, to the best knowledge of the authors, no such minimax estimators have been constructed. In this paper, we provide minimax estimators which perform thresholding simultaneously.

Section 2 develops the new estimator for the canonical form of the model by solving Stein's differential inequality. Sections 3 and 4 provide an approximate Bayesian justification and an empirical Bayes interpretation. Section 5 applies the result to wavelet analysis. The proposed method outperforms several prominent procedures in the statistical wavelet literature. Asymptotic optimality of some of our estimators is established in Section 6.

**2. New estimators for a canonical model.** In this section we shall consider the canonical form of the problem of a multinormal mean estimation problem under squared error loss. Hence we shall assume that our observation

$$Z = (Z_1, \ldots, Z_d) \sim N(\theta, I)$$

has a $d$-dimensional normal distribution with mean $\theta = (\theta_1, \ldots, \theta_d)$, and a known covariance identity matrix $I$. The case when the variance of $Z_i$ is not known will be discussed briefly at the end of Section 5.

The connection of this problem with wavelet analysis will be pointed out in Sections 5 and 6. In short, $Z_i$ and $\theta_i$ represent the wavelet coefficients of the data and the true curve in the same resolution, respectively. Furthermore, $d$ is the dimension of a resolution. For now, we shall seek an estimator of $\theta$ based on $Z$. We shall, without loss of generality, consider an estimator of the form $\delta(Z) = (\delta_1(Z), \ldots, \delta_d(Z))$, where,

$$\delta_i(Z) = Z_i + g_i(Z),$$

where $g(Z): R^d \to R$, and search for $g(Z) = (g_1(Z), \ldots, g_d(Z))$. To insure that the new estimator (perhaps with some thresholding) does better than $Z$ (which does no thresholding), we shall compare the *risk* of $\delta(Z)$ to the risk of $Z$ with respect to the $\ell_2$ norm. Namely,

$$E\|\delta(Z) - \theta\|^2 = E\sum_{i=1}^{d}(\delta_i(Z) - \theta_i)^2.$$

It is obvious that the risk of $Z$ is then $d$. We shall say one is *as good as* the other if the former has a risk no greater than the latter for every $\theta$. Moreover, one *dominates* the other if it is as good as the other and has smaller risk for some $\theta$. Also we shall say that an estimator *strictly dominates* the other if the former has a smaller risk for every $\theta$. Note that $Z$ is a minimax estimator, that is, it minimizes $\sup_\theta E|\delta^0(Z) - \theta|^2$ among all $\delta^0(Z)$. Consequently any $\delta(Z)$ is as good as $Z$ if and only if it is minimax.

To construct an estimator that dominates $Z$, we use the following lemma.

LEMMA 1 [Stein (1981)]. *Suppose that $g: R^d \to R^d$ is a measurable function with $g_i(\cdot)$ as the ith component. If for every $i$, $g_i(\cdot)$ is almost differentiable with respect to the ith component and*

$$E\left(\left|\frac{\partial}{\partial Z_i}g_i(Z)\right|\right) < \infty \quad \text{for } i = 1, \ldots, d,$$

*then*

$$E_\theta\|Z + g(Z) - \theta\|^2 = E_\theta\{d + 2\nabla \cdot g(Z) + \|g(Z)\|^2\},$$

*where $\nabla \cdot g(Z) = \sum_{i=1}^d \frac{\partial g_i(Z)}{\partial Z_i}$. Hence if $g(Z)$ solves the differential inequality*

$$2\nabla \cdot g(Z) + \|g(Z)\|^2 < 0, \tag{1}$$

*the estimator $Z + g(Z)$ strictly dominates $Z$.*

REMARK. $g_i(z)$ is said to be almost differentiable with respect to $z_i$, if for almost all $z_j$, $j \neq i$, $g_i(z)$ can be written as a one-dimensional integral of a function with respect to $z_i$. For such $z_j$'s, $j \neq i$, $g_i(Z)$ is also called absolutely continuous with respect to $z_i$ in Berger (1980).

To motivate the proposed estimator, note that the James–Stein positive estimator has the form

$$\hat{\theta}_i^{\mathrm{JS}} = \left(1 - \frac{d-2}{\|Z\|^2}\right)_+ Z_i$$

with $c_+ = \max(c, 0)$ for any number $c$. This estimator, however, truncates independently of the magnitude of $|Z_i|$. Indeed, it truncates all or none of the coordinates. To construct an estimator that truncates only the coordinates with small $|Z_i|$, it seems necessary to replace $d - 2$ by a decreasing function $h(|Z_i|)$ of $|Z_i|$ and consider

$$\hat{\theta}_i^+ = \left(1 - \frac{h(|Z_i|)}{D}\right)_+ Z_i,$$

where $D$, independent of $i$, is yet to be determined. [In a somewhat different approach, Beran and Dümbgen (1998) construct a modulation estimator corresponding to a monotonic shrinkage factor.] With such a form, $\hat{\theta}_i^+ = 0$ if $h(|Z_i|) \geq D$, which has a better chance of being satisfied when $|Z_i|$ is small.

We consider the simple choice $h(|Z_i|) = a|Z_i|^{-2/3}$, and let $D = \Sigma |Z_i|^{4/3}$. This leads to the untruncated version $\hat{\theta}$ with the $i$th component

$$(2) \qquad \hat{\theta}_i(Z) = Z_i + g_i(Z) \qquad \text{where } g_i(Z) = -aD^{-1}\text{sign}(Z_i)|Z_i|^{1/3}.$$

Here and later $\text{sign}(Z_i)$ denotes the sign of $Z_i$. It is possible to use other decreasing functions $h(|Z_i|)$ and other $D$.

In general, we consider, for a fixed $\beta \leq 2$, an estimator of the form

$$(3) \qquad \hat{\theta}_i = Z_i + g_i(Z),$$

where

$$(4) \qquad g_i(Z) = -a\frac{\text{sign}(Z_i)|Z_i|^{\beta-1}}{D} \quad \text{and} \quad D = \sum_{i=1}^d |Z_i|^\beta.$$

Although at first glance it may seem hard to justify this estimator, it has a Bayesian and empirical Bayes justification in Sections 3 and 4. It contains, as a special case with $\beta = 2$, the James–Stein estimator. Now we have:

THEOREM 2.  *For $d \geq 3$ and $1 < \beta \leq 2$, $\hat{\theta}(Z)$ is minimax if and only if*

$$0 < a \leq 2(\beta - 1)\inf_\theta \frac{E_\theta(D^{-1}\sum_{i=1}^d |Z_i|^{\beta-2})}{E_\theta(D^{-2}\sum_{i=1}^d |Z_i|^{(2\beta-2)})} - 2\beta.$$

PROOF.  Obviously for $Z_j \neq 0$, $\forall j \neq i$, $g_i(Z)$ can be written as the one-dimensional integral of

$$\frac{\partial}{\partial Z_i} g_i(Z) = \beta(-a)(-1)D^{-2}|Z_i|^{(2\beta-2)} + (\beta-1)(-a)D^{-1}(|Z_i|^{\beta-2})$$

with respect to $Z_i$. (The only concern is at $Z_i = 0$.) Consider only nonzero $Z_j$'s, $j \neq i$. Since $\beta > 1$, this function, however, is integrable with respect to $Z_i$ even over an integral including zero. It takes some effort to prove that $E(|\frac{\partial}{\partial Z_i}g_i(Z)|)$ is finite. However, one only needs to focus on $Z_j$ close to zero. Using the spherical-like transformation $r^2 = \sum |Z_i|^\beta$, we may show that if $d \geq 3$ and $\beta > 1$, both terms in the above displayed expression are integrable.

Now
$$\|g(Z)\|^2 = a^2 D^{-2} \sum_{i=1}^{d} |Z_i|^{2\beta-2}.$$

Hence
$$E_\theta \|Z + g(Z) - \theta\|^2 \leq d \quad \text{for every } \theta$$

if and only if
$$E_\theta \{2\nabla \cdot g(Z) + \|g(Z)\|^2\} \leq 0 \quad \text{for every } \theta,$$

that is,

(5) $$E_\theta \bigg( a \bigg( (2\beta) D^{-2} \sum_{i=1}^{d} |Z_i|^{(2\beta-2)} - (2\beta-2) D^{-1} \sum_{i=1}^{d} |Z_i|^{\beta-2} \bigg) + a^2 D^{-2} \sum_{i=1}^{d} |Z_i|^{2\beta-2} \bigg) \leq 0$$
$$\text{for every } \theta,$$

which is equivalent to the condition stated in the theorem. □

THEOREM 3. *The estimator $\hat{\theta}(Z)$ with the ith component given in (3) and (4) is minimax provided $0 < a \leq 2(\beta-1)d - 2\beta$ and $1 < \beta \leq 2$. Unless $\beta = 2$ and $a$ is taken to the upper bound, $\hat{\theta}(Z)$ dominates $Z$.*

PROOF. By the correlation inequality
$$d\bigg(\sum_{i=1}^{d} |Z_i|^{2\beta-2}\bigg) \leq \bigg(\sum_{i=1}^{d} |Z_i|^{\beta-2}\bigg)\bigg(\sum_{i=1}^{d} |Z_i|^\beta\bigg).$$

Strict inequality holds almost surely if $\beta < 2$. Hence
$$\frac{E_\theta(D^{-1} \sum_{i=1}^{d} |Z_i|^{\beta-2})}{E_\theta(D^{-2} \sum_{i=1}^{d} |Z_i|^{2\beta-2})} \geq \frac{E_\theta D^{-1} \sum |Z_j|^{\beta-2}}{(1/d) E_\theta D^{-1} \sum |Z_i|^{\beta-2}} = d.$$

Hence if $0 < a \leq 2(\beta-1)d - 2\beta$, the condition in Theorem 2 is satisfied, implying minimaxity of $\hat{\theta}(Z)$. The rest of the statement of the theorem is now obvious. □

The following theorem is a generalization of Theorem 6.2 on page 302 of Lehmann (1983) and Theorem 5.4 on page 356 of Lehmann and Casella

(1998). It shows that taking the positive part will improve the estimator componentwise. Specifically for an estimator $(\tilde{\theta}_1(Z), \ldots, \tilde{\theta}_d(Z))$ where

$$\tilde{\theta}_i(Z) = (1 - h_i(Z))Z_i,$$

the positive part estimator of $\tilde{\theta}_i(Z)$ is denoted as

$$\tilde{\theta}_i^+(Z) = (1 - h_i(Z))_+ Z_i.$$

THEOREM 4. *Assume that $h_i(Z)$ is symmetric with respect to the ith coordinate. Then*

$$E_\theta(\theta_i - \tilde{\theta}_i^+)^2 \leq E_\theta(\theta_i - \tilde{\theta}_i)^2.$$

*Furthermore, if*

(6) $$P_\theta(h_i(Z) > 1) > 0,$$

*then*

$$E_\theta(\theta_i - \tilde{\theta}_i^+)^2 < E_\theta(\theta_i - \tilde{\theta}_i)^2.$$

PROOF. Simple calculation shows that

(7) $$E_\theta(\theta_i - \tilde{\theta}_i^+)^2 - E_\theta(\theta_i - \tilde{\theta}_i)^2 = E_\theta((\tilde{\theta}_i^+)^2 - \tilde{\theta}_i^2) - 2\theta_i E_\theta(\tilde{\theta}_i^+ - \tilde{\theta}_i).$$

Let us calculate the expectation by conditioning on $h_i(Z)$. For $h_i(Z) \leq 1$, $\tilde{\theta}_i^+ = \tilde{\theta}_i$. Hence it is sufficient to condition on $h_i(z) = b$ where $b > 1$ and show that

$$E_\theta((\tilde{\theta}_i^+)^2 - \tilde{\theta}_i^2 | h_i(Z) = b) - 2\theta_i E_\theta(\tilde{\theta}_i^+ - \tilde{\theta}_i | h_i(Z) = b) \leq 0,$$

or equivalently,

$$-E_\theta(\tilde{\theta}_i^2 | h_i(Z) = b) + 2\theta_i(1 - b)E_\theta(Z_i | h_i(Z) = b) \leq 0.$$

Obviously, the last inequality is satisfied if we can show

$$\theta_i E_\theta(Z_i | h_i(Z) = b) \geq 0.$$

We may further condition on $Z_j = z_j$ for $j \neq i$ and it suffices to establish

(8) $$\theta_i E_\theta(Z_i | h_i(Z) = b, Z_j = z_j, j \neq i) \geq 0.$$

Given that $Z_i = z_j$, $j \neq i$, consider only the case where $h_i(Z) = b$ has solutions. Due to symmetry of $h_i(Z)$, these solutions are in pairs. Let $\pm y_k$, $k \in K$, denote the solutions. Hence the left-hand side of (8) equals

$$\theta_i E_\theta(Z_i | Z_i = \pm y_k, k \in K)$$
$$= \sum_{k \in K} \theta_i E_\theta(Z_i | Z_i = \pm y_k) P_\theta(Z_i = \pm y_k | Z_i = \pm y_k, k \in K).$$

Note that

$$(9) \quad \theta_i E_\theta(Z_i | Z_i = \pm y_k) = \frac{\theta_i y_k e^{y_k \theta_i} - \theta_i y_k e^{-y_k \theta_i}}{e^{y_k \theta_i} + e^{-y_k \theta_i}},$$

which is symmetric in $\theta_i y_k$ and is increasing for $\theta_i y_k > 0$. Hence (9) is bounded below by zero, a bound obtained by substituting $\theta_i y_k = 0$ in (9). Consequently we establish that (7) is nonpositive, implying that $\tilde{\theta}^+$ in as good as $\tilde{\theta}$.

The strict inequality of the theorem can be established by noting that the right-hand side of (7) is bounded above by $E_\theta[(\tilde{\theta}_i^+)^2 - \tilde{\theta}_i^2]$ which by (6) is strictly negative. □

Theorem 4 implies the following corollary.

COROLLARY 5. *Under the assumptions on $a$ and $\beta$ in Theorem 3, $\hat{\theta}^+$ with the ith component*

$$(10) \quad \hat{\theta}_i^+ = (1 - aD^{-1}|Z_i|^{\beta-2})_+ Z_i$$

*strictly dominates $Z$.*

It is interesting to note that the estimator (10), for $\beta < 2$, does give zero as the estimator when the $|Z_i|$ are small. When applied to wavelet analysis, it truncates the small wavelet coefficients and shrinks the large wavelet coefficients. The estimator lies in a data-chosen reduced model.

Moreover, for $\beta = 2$, Theorem 3 reduces to the classical result of Stein (1981) and (10) to the positive part James–Stein estimator. The upper bound of $a$ for domination stated in Theorem 3 works only if $\beta > 1$ and $d > \beta/(\beta - 1)$. We know that for $\beta \leq \frac{1}{2}$, $\hat{\theta}$ fails to dominate $Z$ because of the calculations leading to (11) below. We are unable to prove that $\hat{\theta}$ dominates $Z$ for $\frac{1}{2} < \beta \leq 1$. However, for such $\beta$'s, $\hat{\theta}$ has a smaller Bayes risk than $Z$ if the condition (11) below is satisfied.

*A remark about an explicit formula for $a$.* In wavelet analysis, a vast majority of the wavelet coefficients of a reasonably smooth function are zero. Consequently, it seems good to choose an estimator that shrinks a lot and hence using $a$ larger than the upper bound in Theorem 3 is desirable. Although Theorem 2 provides the largest possible $a$ for domination in the frequentist sense, the bound is difficult to evaluate in computation and hence difficult to use in a real application. Hence we took an alternative approach by assuming that $\theta_i$ are independently identically distributed (i.i.d.) $N(0, \tau^2)$. It can be shown by a detailed calculation [see Zhou and Hwang

(2003)] that the estimator (3) and (4) has a smaller Bayes risk than $Z$ for all $\tau^2$ if and only if

$$
(11) \qquad 0 < a < a_\beta = 2 \Big/ E\left[\sum_{i=1}^{d} |\xi_i|^{2\beta-2} \Big/ \left(\sum_{i=1}^{d} |\xi_i|^\beta\right)^2\right],
$$

where $\xi_i$ are i.i.d. standard normal random variables.

What is the value of $a_\beta$? It is easy to numerically calculate the bound $a_\beta$ by simulating $\xi_i$, which we did for $a$ up to 100. It is shown that $a_\beta$, $\beta = \frac{4}{3}$, is at least as big as $(5/3)(d-2)$. Using Berger's (1976) tail minimaxity argument, we come to the conclusion that $\hat{\theta}^+$, with the $i$th component

$$
(12) \qquad \hat{\theta}_i = \left(1 - \frac{(5/3)(d-2)Z_i^{-2/3}}{\sum_{i=1}^d Z_i^{4/3}}\right)_+ Z_i,
$$

would possibly dominate $Z$. For various $d$'s including $d = 50$ this was shown to be true numerically.

To derive a general formula for $a_\beta$ for all $\beta$, we then establish that the limit of $a_\beta/d$ as $d \to \infty$ equals, for $1/2 < \beta < 2$,

$$
(13) \qquad C_\beta = 4[\Gamma((\beta+1)/2)]^2 / [\sqrt{\pi}\Gamma((2\beta-1)/2)].
$$

It may be tempting to use $(d-2)C_\beta$. However, we recommend

$$
(14) \qquad a = (0.97)(d-2)C_\beta,
$$

so that at $\beta = 4/3$, (14) becomes $(5/3)(d-2)$. Berger's tail minimaxity argument and many numerical studies indicate that this $a$ enables (10) to have a better risk than $Z$.

**3. Approximate Bayesian justification.** It would seem interesting to justify the proposed estimation from a Bayesian point of view. To do so, we consider a prior of the form

$$
\pi(\theta) = \begin{cases} 1, & \|\theta\|_\beta \leq 1, \\ 1/(\|\theta\|_\beta)^{\beta c}, & \|\theta\|_\beta > 1, \end{cases}
$$

where $\|\theta\|_\beta = (\sum \|\theta_i\|^\beta)^{1/\beta}$, and $c$ is a positive constant which can be specified to match the constant $a$ in (10). In general the Bayes estimator is given by

$$
Z + \nabla \log m(Z),
$$

where $m(Z)$ is the marginal probability density function of $Z$, namely,

$$
m(Z) = \int \cdots \int \frac{e^{-(1/2)\|Z-\theta\|^2}}{(\sqrt{2\pi})^d} \pi(\theta)\,d\theta.
$$

The following approximation follows from Brown (1971), which asserts that $\nabla \log m(Z)$ can be approximated by $\nabla \log \pi(Z)$. The proof is given in the Appendix.

THEOREM 6. *With $\pi(\theta)$ and $m(X)$ given above,*

$$\lim_{|Z_i|\to+\infty} \frac{\nabla_i \log m(Z)}{\nabla_i \log \pi(Z)} = 1.$$

Hence by Theorem 6, the $i$th component of the Bayes estimator equals approximately

$$Z_i + \nabla_i \log \pi(Z) = Z_i - \frac{c\beta|Z_i|^{\beta-1}\operatorname{sign}(Z_i)}{\sum |Z_i|^\beta}.$$

This is similar to the untruncated version of $\hat{\theta}$ in (2) and (3).

**4. Empirical Bayes justification.** Based on several signals and images, Mallat (1989) proposed a prior for the wavelet coefficients $\theta_i$ as the exponential power distribution with the probability density function (p.d.f.) of the form

(15) $$f(\theta_i) = ke^{-|\theta_i/\alpha|^\beta},$$

where $\alpha$ and $\beta < 2$ are positive constants and

$$k = \beta/(2\alpha\Gamma(1/\beta))$$

is the normalization constant. See also Vidakovic [(1999), page 194]. Using the method of moments, Mallat estimated the values of $\alpha$ and $\beta$ to be 1.39 and 1.14 for a particular graph. However, $\alpha$ and $\beta$ are typically unknown.

It seems reasonable to derive an empirical Bayes estimator based on this class of prior distributions. First we assume that $\alpha$ is known. Then the Bayes estimator of $\theta_i$ is

$$Z_i + \frac{\partial}{\partial Z_i}\log m(Z).$$

Similar to the argument in Theorem 6 and noting that for $\beta < 2$,

$$e^{-|\theta_i+Z_i|^\beta/\alpha^\beta}/e^{-|\theta_i|^\beta/\alpha^\beta} \to 1 \qquad \text{as } \theta_i \to \infty,$$

the Bayes estimator can be approximated by

(16) $$Z_i + \frac{\partial}{\partial Z_i}\log \pi(Z_i) = Z_i - \frac{\beta}{\alpha^\beta}|Z_i|^{\beta-1}\operatorname{sign}(Z_i).$$

Note that, under the assumption that $\alpha$ is known, the above expression is also the asymptotic expression of the maximum likelihood estimator of $\theta_i$ by maximizing the joint p.d.f. of $(Z_i, \theta_i)$. See Proposition 1 of Antoniadis, Leporini and Pesquet (2002) as well as (8.23) of Vidakovic (1999). In the latter reference, the sign of $Z_i$ of (16) is missing due to a minor typographical error.

Since $\alpha$ is unknown, it seems reasonable to replace $\alpha$ in (16) by an estimator. Assume that $\theta_i$'s are observable. Then by (15) the joint density of $(\theta_1, \ldots, \theta_d)$ is

$$\left[\frac{\beta}{2\alpha\Gamma(1/\beta)}\right]^d e^{-\Sigma(|\theta_i|^\beta/\alpha^\beta)}.$$

Differentiating this p.d.f. with respect to $\alpha$ gives the maximum likelihood estimator of $\alpha^\beta$ as

(17) $$(\beta\Sigma|\theta_i|^\beta)/d.$$

However, $\theta_i$ is unknown and hence the above expression can be further estimated by (16). For $\beta < 2$, the second term in (16) has a smaller order than the first when $|Z_i|$ is large. Replacing $\theta_i$ by the dominating first term $Z_i$ in (16) leads to an estimator of $\alpha^\beta$ as $(\beta\Sigma|Z_i|^\beta)/d$.

Substituting this into (16) gives

$$Z_i - \frac{d}{\Sigma|Z_i|^\beta}|Z_i|^{\beta-1}\operatorname{sign}(Z_i),$$

which is exactly the estimator $\hat{\theta}_i$ in (2) and (3) with $a = d$. Hence we have succeeded in deriving $\hat{\theta}_i$ as an empirical Bayes estimator when $Z_i$ is large.

**5. Connection to the wavelet analysis and the numerical results.** Wavelets have become a very important tool in many areas including mathematics, applied mathematics, statistics and signal processing. They are also applied to numerous other areas of science such as chemometrics and genetics.

In statistics, wavelets have been applied to function estimation with amazing results of being able to catch the sharp change of a function. Celebrated contributions by Donoho and Johnstone (1994, 1995) focus on developing thresholding techniques and asymptotic theories. In the 1994 paper, relative to the oracle risk, their VisuShrink is shown to be asymptotically optimal. Further in the 1995 paper, the expected squared error loss of their SureShrink is shown to achieves the global asymptotic minimax rate over Besov spaces. Cai (1999) improved on their result by establishing that Block James–Stein (BlockJS) thresholding achieves exactly the asymptotic global or local minimax rate over various classes of Besov spaces.

Now specifically let $Y = (Y_1, \ldots, Y_n)'$ be samples of a function $f$, satisfying

(18) $$Y_i = f(t_i) + \varepsilon_i,$$

where $t_i = (i-1)/n$ and $\varepsilon_i$ are i.i.d. $N(0, \sigma^2)$. Here $\sigma^2$ is assumed to be known and is taken to be 1 without loss of generality. See a comment at the end of the paper regarding the unknown $\sigma$ case. One wishes to choose an estimate $\hat{f} = (\hat{f}(t_1), \ldots, \hat{f}(t_n))$ so that its risk function

$$E\|\hat{f} - f\|^2 = E \sum_{i=1}^{n} (\hat{f}(t_i) - f(t_i))^2 \tag{19}$$

is as small as possible. Many discrete wavelet transformations are orthogonal transformations. See Donoho and Johnstone (1995). Consequently, there exists an orthogonal matrix $W$ such that the wavelet coefficients of $Y$ and $f$ are $Z = WY$ and $\theta = Wf$. Obviously the components $Z_i$ of $Z$ are independent, having a normal distribution with mean $\theta_i$ and standard deviation 1. Hence previous sections apply and exhibit many good estimators of $\theta$. Note that, by orthogonality of $W$, for any estimator $\delta(Z)$ of $\theta$, its risk function is identical to $W'\delta(Z)$ as an estimator of $f = W'\theta$. Hence the good estimators in previous sections can be inversely transformed to estimate $f$ well.

In all the applications to wavelets discussed in this paper, the estimators (including our proposed estimator) apply separately to the wavelet coefficients of the same resolution. Hence in (12), for example, $d$ is taken to be the number of coefficients of a resolution when applied to the resolution. In all the literature that we are aware of, this has been the case as well.

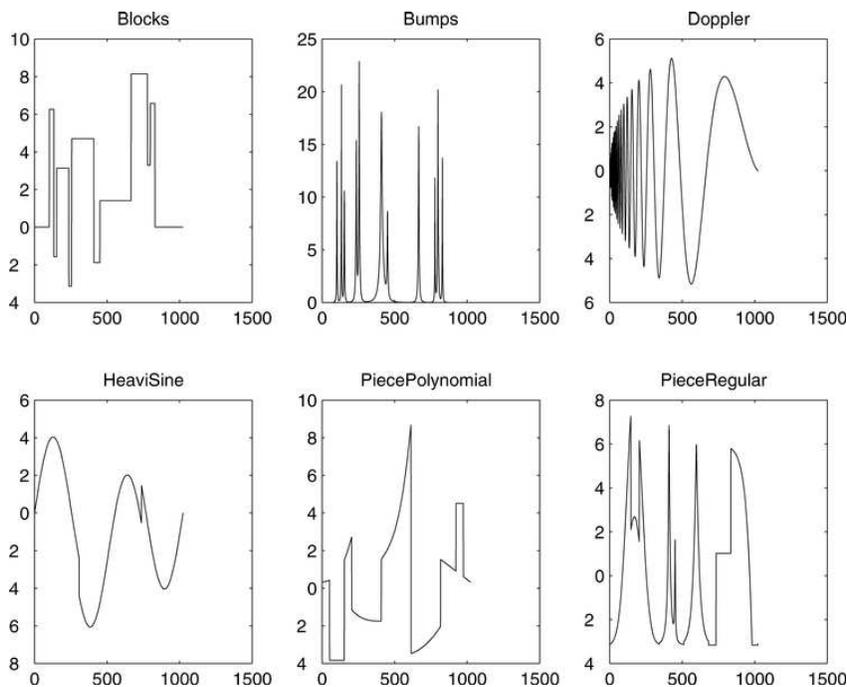

FIG. 1. *The curves represent the true curves $f(t)$ in (18).*

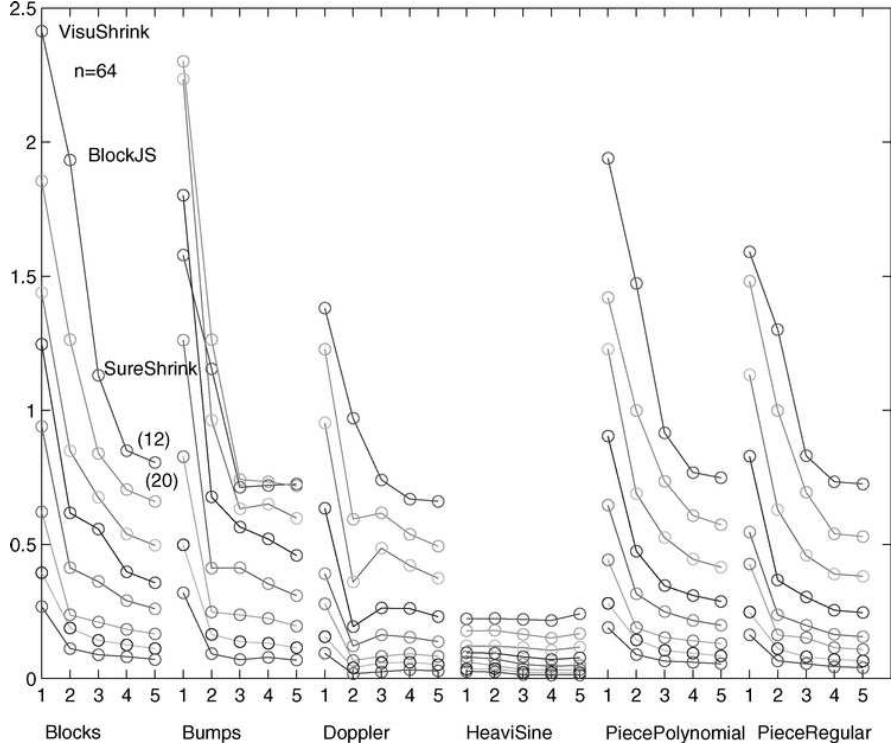

FIG. 2. *In each of the six cases corresponding to Blocks, Bumps, and so on, the eight curves plot the risk functions, from top to bottom, when $n = 64, 128, \ldots, 8192$. For each curve (see, e.g., the top curve on the left), the circles "o" from left to right give, with respect to $n$, the relative risks of VisuShrink, Block James–Stein, SureShrink and the proposed methods* (12) *and* (20).

In addition to considering the estimator (12), which is a special case of (10) with $\beta = 4/3$, we also propose a modification (10) with an estimated $\beta$. The estimator $\hat{\beta}$ for $\beta$ is constructed by minimizing, for each resolution, the Stein unbiased risk estimator (SURE) for the risk of (10). The quantity SURE is basically the expression inside the expectation on the right-hand side of (A.4) summing over $i$, $1 \leq i \leq d$, except that $a$ is replaced by $a_\beta$. [Note that $D$ in (A.4) depends on $\beta$ as well.] The resultant estimator is denoted as

(20) $\qquad\qquad\qquad \hat{\theta}^S = (10) \qquad \text{with } \beta \text{ replaced by } \hat{\beta}.$

Figure 1 gives six true curves (made famous by Donoho and Johnstone) from which the data are generated. For these six cases, Figure 2 plots the ratios of the risks of the aforementioned estimators to $n$, the risk of $Y$. Since most relative risks are less than 1, this indicates that most estimators perform better than the raw data $Y$. Our estimators $\hat{\theta}^+$ in (12) and $\hat{\theta}^S$ in (20), however, are the ones that are consistently better than $Y$. Furthermore, our estimators $\hat{\theta}^+$ and $\hat{\theta}^S$ virtually dominate all the other estimators in risk. Generally, $\hat{\theta}^S$ performs better than $\hat{\theta}^+$ virtually in all cases.

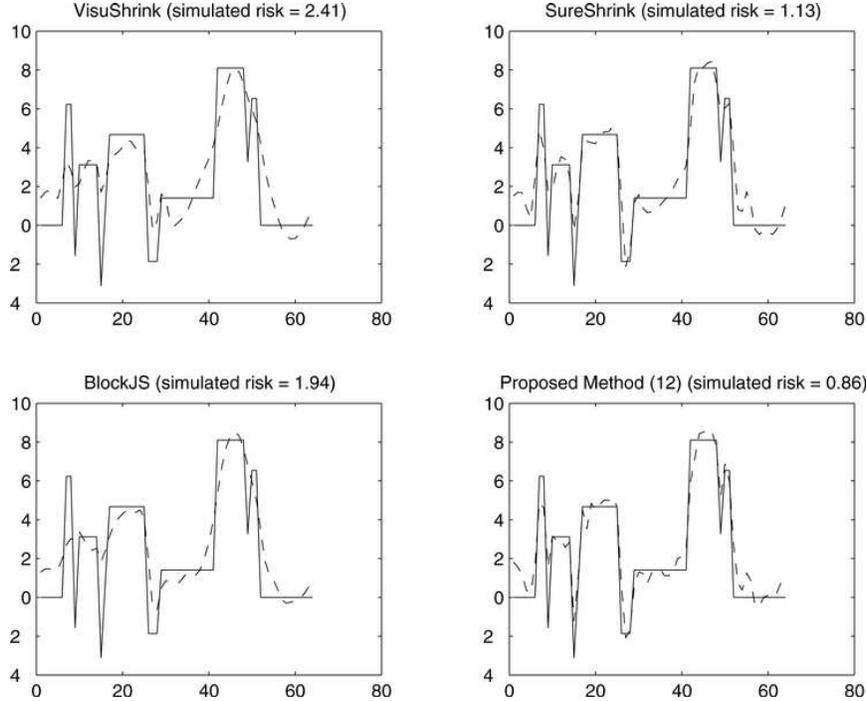

FIG. 3. *Solid lines represent the true curves, and dotted lines represent the curves corresponding to various estimators. The simulated risk is based on* 500 *simulations.*

As shown in Figure 2, the difference in risks between $\hat{\theta}^+$ and $\hat{\theta}^S$ is quite minor. Since $\hat{\theta}^+$ is computationally less intensive, we focus on $\hat{\theta}^+$ for the rest of the numerical studies.

Picturewise, our estimator does slightly better than other estimators. See Figure 3 for an example. Note that the picture corresponding to $\hat{\theta}^+$ distinguishes most clearly the first and second bumps from the right.

Based on asymptotic calculation, the next section also recommends a choice of $a$ in (21). It would seem interesting to comment on its numerical performance. The difference between the $a$'s defined in (14) and (22) is very small when $64 \leq n \leq 8192$ and when $\beta$ is estimated by minimizing SURE. Consequently, for such $\beta$, the risk functions of the two estimators with different $a$'s are very similar, with a difference virtually bounded by 0.02. The finite sample estimator [where $a$ is defined in (14)] has a smaller risk about 75% of the time.

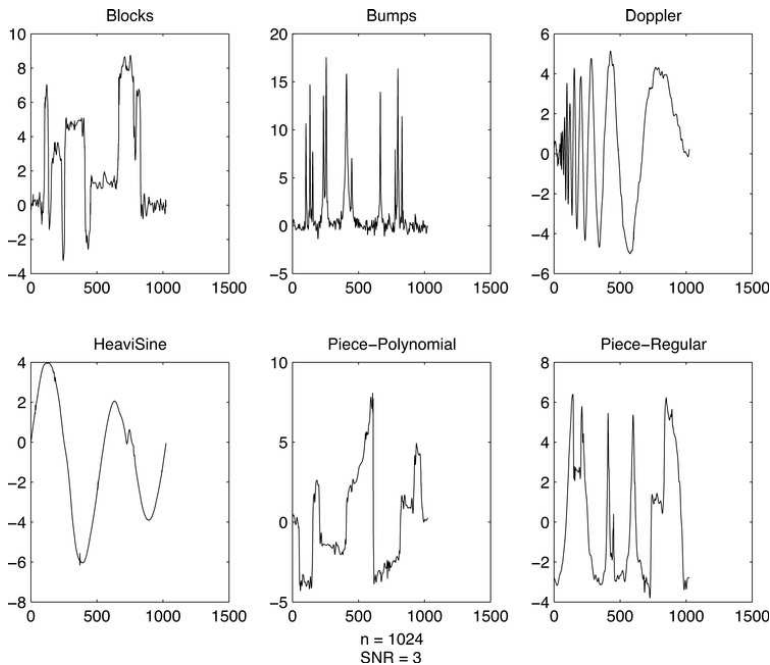

Fig. 4.  *Proposed estimator* (12) *applied to reconstruct Figure* 1.

The James–Stein estimator produces very attractive risk functions, sometimes as good as the proposed estimator (12). However, it does not seem to produce good graphs. Compare Figures 4 and 5.

In the simulation studies we use the procedures MultiVisu and MultiHybrid which are VisuShrink and SureShrink in WaveLab802. See http://www-stat. stanford.edu/∼wavelab. We use Symmlet 8 to do wavelet transformation. In Figure 2 signal-to-noise ratio (SNR) is taken to be 3. Results are similar for other SNRs. To include the block thresholding result of Cai (1999), we choose the lowest integer resolution level $j \geq \log_2(\log n) + 1$.

**A comment about the case where $\sigma^2$ is not known to be 1.**  When $\sigma$ is known and is not equal to 1, a simple transformation applied to the problem suggests that (10) be modified with $a$ replaced by $a\sigma^2$. When $\sigma$ is unknown, one could then estimate $\sigma$ by $\hat{\sigma}$, the proposed estimator for $\sigma$ in Donoho and Johnstone [(1995), page 1218]. With this modification in (12) (and even with the SURE estimated $\beta$), the resultant estimators are not minimax according to some numerical simulations. However, they still perform the best or nearly the best among all the estimators studied in Figure 2.

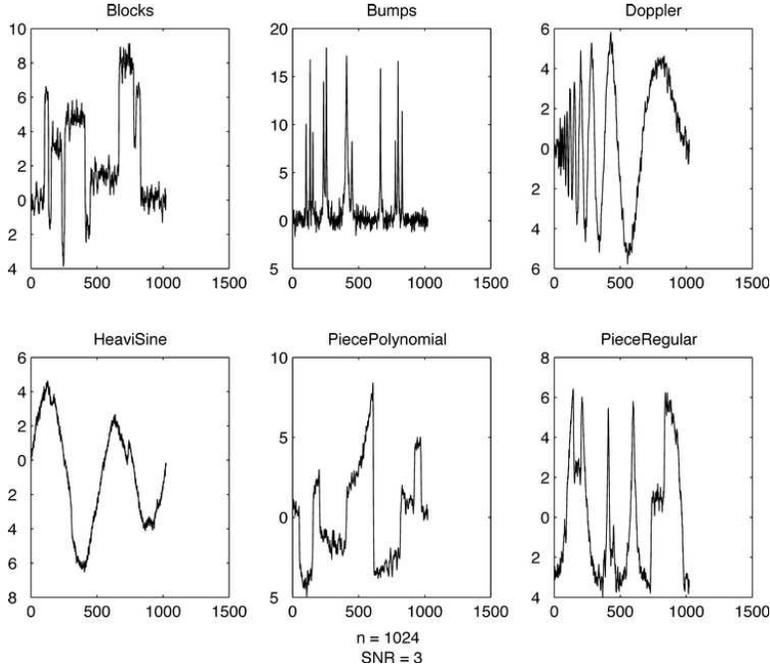

FIG. 5. *James–Stein positive part applied to reconstruct Figure* 1.

**6. Asymptotic optimality.** To study the asymptotic rate of a wavelet analysis estimator, it is customary to assume the model

$$Y_i = f(t_i) + \varepsilon_i, \qquad i = 1, \ldots, n, \tag{21}$$

where $t_i = (i-1)/n$ and $\varepsilon_i$ are assumed to be i.i.d. $N(0,1)$. The estimator $\hat{f}$ for $f(\cdot)$ that can be proved asymptotically optimal applies estimator (10) with

$$a = d(2\ln d)^{(2-\beta)/2} m_\beta, \qquad 0 \leq \beta \leq 2, \tag{22}$$

and

$$m_\beta = E|\varepsilon_i|^\beta = 2^{\beta/2}\Gamma((\beta+2)/2)\sqrt{\pi}$$

to the wavelet coefficients $Z_i$ of each resolution with dimensionality $d$ of the wavelet transformation of the $Y_i$'s. After applying the estimator to each resolution one at a time to come up with the new wavelet coefficient estimators, one then uses the wavelet base function to obtain one function $\hat{f}$ in the usual way.

To state the theorem, we use $B_{p,q}^\alpha$ to denote the Besov space with smoothness $\alpha$ and shape parameters $p$ and $q$. The definition of the Besov class $B_{p,q}^\alpha(M)$ with respect to the wavelet coefficients is given in (A.19). Now the asymptotic theorem is given below.

THEOREM 7. *Assume that the wavelet $\psi$ is t-regular, that is, $\psi$ has t vanishing moments and t continuous derivatives. Then there exists a constant $C$ independent of $n$ and $f$ such that*

$$(23) \quad \sup_{f \in B^{\alpha}_{p,q}(M)} E \int_0^1 |f(t) - \hat{f}(t)|^2 \, dt \leq C(\ln n)^{1-\beta/2} n^{-2\alpha/(2\alpha+1)},$$

*for all $M > 0$, $0 < \alpha < r$, $q \geq 1$ and $p > \max(\beta, \frac{1}{\alpha}, 1)$.*

The asymptotic optimality stated in (23) is as good as what has been established for hard and soft thresholding estimators in Donoho and Johnstone (1994), the Garrott method in Gao (1998) and Theorem 4 in Cai (1999) and the SCAD method in Antoniadis and Fan (2001). However, the real advantage of our estimator is in the finite sample risk as reported in Section 5. Also our estimators are constructed to be minimax and hence have finite risk functions uniformly smaller than the risk of $Z$. This estimator $\hat{\theta}^A$ for $\beta = 4/3$, however, has a risk very similar to (12). See Section 5.

## APPENDIX

PROOF OF THEOREM 6. Assume that $|Z_i| > 1$. We have

$$\lim_{|Z_i| \to \infty} \frac{\nabla_i \log m(Z)}{\nabla_i \log \pi(Z)} = \lim_{|Z_i| \to +\infty} \frac{\pi(Z)}{m(Z)} \cdot \frac{(\partial/\partial Z_i) m(Z)}{(\partial/\partial Z_i) \pi(Z)}.$$

We shall prove only

$$\lim_{|Z_i| \to \infty} \frac{m(Z)}{\pi(Z)} = 1,$$

since

$$\lim_{|Z_i| \to \infty} \frac{(\partial/\partial Z_i) m(Z)}{(\partial/\partial Z_i) \pi(Z)} = 1$$

can be similarly established.

Now

$$m(Z) = \int \cdots \int \frac{1}{(\sqrt{2\pi})^p} e^{-(1/2)\|Z-\theta\|^2} \pi(\theta) \, d\theta$$

$$= \int \cdots \int_{\|\theta\|_\beta \leq 1} \frac{1}{(\sqrt{2\pi})^p} e^{-(1/2)\|Z-\theta\|^2} \, d\theta$$

$$\quad + \int \cdots \int_{\|\theta\|_\beta > 1} \frac{1}{(\sqrt{2\pi})^p} e^{-(1/2)\|Z-\theta\|^2} \frac{1}{\|\theta\|_\beta^{\beta c}} \, d\theta$$

$$= m_1 + m_2, \quad \text{say.}$$

Obviously, as $|Z_i| \to +\infty$, $m_1$ has an exponentially decreasing tail. Hence

$$\lim_{|Z_i| \to +\infty} \frac{m_1}{\pi(Z)} = 0.$$

By the change of variable $\theta = Z + y$, we have

$$m_2/\pi(Z) = \int \cdots \int_{\|Z+y\|_\beta > 1} \frac{1}{(\sqrt{2\pi})^p} e^{-(1/2)\|y\|^2} \frac{\|Z\|_\beta^{\beta c}}{\|Z+y\|_\beta^{\beta c}} \, dy.$$

To prove the theorem, it suffices to show the above expression converges to 1. In doing so, we shall apply the dominated convergence theorem to show that we may pass to the limit inside the above integral. After passing to the limit, it is obvious that the integral becomes 1.

The only argument left is to show that the dominated convergence theorem can be applied. To do so, we seek an upper bound $F(y)$ for

$$\|Z\|_\beta^{\beta c} / \|Z+y\|_\beta^{\beta c} \quad \text{when } \|Z+y\|_\beta > 1.$$

Now for $\|Z+y\|_\beta > 1$,

$$\|Z\|_\beta^{\beta c} \leq C_p(\|Z+y\|_\beta^{\beta c} + \|y\|_\beta^{\beta c}),$$

that is,

$$\frac{\|Z\|_\beta^{\beta c}}{\|Z+y\|_\beta^{\beta c}} \leq C_p\left(1 + \frac{\|y\|_\beta^{\beta c}}{\|Z+y\|_\beta^{\beta c}}\right) \leq C_p(1 + \|y\|_\beta^{\beta c}).$$

Hence if we take $C_p(1 + \|y\|_\beta^{\beta c})$ as $F(y)$, then

$$\int \cdots \int_{\|Z+y\|_\beta > 1} \frac{1}{(\sqrt{2\pi})^p} e^{-(1/2)\|y\|^2} F(y) \, dy < +\infty.$$

Consequently, we may apply the dominated convergence theorem, which completes the proof. □

PROOF OF THEOREM 7. Before relating to model (21), we shall work on the canonical form:

$$Z_i = \theta_i + \sigma \varepsilon_i, \qquad i = 1, 2, \ldots, d,$$

where $\sigma > 0$ and the $\varepsilon_i$'s are independently identically distributed standard normal random errors. Here $\hat{\theta} = (\hat{\theta}_1, \ldots, \hat{\theta}_d)$ denotes the estimator in (10) with $a$ defined in (22). For the rest of the paper $C$ denotes a generic quantity independent of $d$ and the unknown parameters. Hence the $C$'s below are not necessarily identical.

We shall first prove Lemma A.1 below. Inequality (A.1) will be applied to the lower resolutions in the wavelet regression. The other two inequalities (A.2) and (A.3) are for higher resolutions. □

LEMMA A.1. *For any $0 < \beta < 2$, $0 < \delta < 1$, and some $C > 0$, independent of $d$ and the $\theta_i$'s, we have*

$$\sum_{i=1}^{d} E(\hat{\theta}_i - \theta_i)^2 \leq C\sigma^2 d (\ln d)^{(2-\beta)/2} \tag{A.1}$$

*and*

$$E(\hat{\theta}_i - \theta_i)^2 \leq C(\theta_i^2 + \sigma^2 d^{\delta-1}(\ln d)^{-1/2}) \tag{A.2}$$

$$\text{if } \sum_{i=1}^{d} |\theta_i|^\beta \leq \sigma^\beta \left(\frac{2-\beta}{2\beta}\right)^\beta \delta^2 m_\beta d.$$

Here and below $m_\beta$ denotes the expectation of $|\varepsilon_i|^\beta$, defined right above the statement of Theorem 7. Furthermore, for any $0 < \beta < 1$, there exists $C > 0$ such that

$$E(\hat{\theta}_i - \theta_i)^2 \leq C\sigma^2 \ln d. \tag{A.3}$$

PROOF OF LEMMA A.1. Without loss of generality we will prove the theorem for the case $\sigma = 1$. By Stein's identity,

$$E(\hat{\theta}_i - \theta_i)^2 = E\bigg[1 + (Z_i^2 - 2)I_i \tag{A.4}$$

$$+ \left(\frac{a^2 |Z_i|^{2\beta-2}}{D^2} - 2a(\beta-1)\frac{|Z_i|^{\beta-2}}{D} + 2a\beta\frac{|Z_i|^{2\beta-2}}{D^2}\right)I_i^c\bigg].$$

Here $I_i$ denotes the indicator function $I(a|Z_i|^{\beta-2} > D)$ and $I_i^c = 1 - I_i$. Consequently

$$I_i = 1 \quad \text{if } |Z_i|^{2-\beta} < a/D \tag{A.5}$$

and

$$I_i^c = 1 \quad \text{if } a|Z_i|^{\beta-2}/D \leq 1. \tag{A.6}$$

From (A.4) and after some straightforward calculations,

$$E\sum_{i=1}^{d}(\hat{\theta}_i - \theta_i)^2$$

$$= d + E\bigg[\sum_{i=1}^{d}(|Z_i|^{2-\beta}|Z_i|^\beta - 2)I_i \tag{A.7}$$

$$+ \frac{a|Z_i|^{\beta-2}}{D}\left(\frac{a|Z_i|^\beta}{D} - 2(\beta-1) + 2\beta\frac{|Z_i|^\beta}{D}\right)I_i^c\bigg].$$

Using this and the upper bounds in (A.5) and (A.6), we conclude that (A.7) is bounded above by

$$d + E\left[\sum_{i=1}^{d} \frac{a|Z_i|^\beta}{D} + \frac{a|Z_i|^\beta}{D} + 2\beta\frac{|Z_i|^\beta}{D}\right] + 2|\beta - 1|d \leq C(\ln d)^{(2-\beta)/2}d,$$

completing the proof of (A.1).

To derive (A.2), note that

$$E(1 + (Z_i^2 - 2)I_i) = \theta_i^2 + E(-Z_i^2 + 2)I_i^c.$$

This and (A.4) imply that

$$E(\hat{\theta}_i - \theta_i)^2 = \theta_i^2 + E\left\{\left[\left(\frac{a|Z_i|^{\beta-2}}{D}\right)^2 Z_i^2 - Z_i^2\right]I_i^c\right\}$$

$$+ E\left\{\left[-2(\beta-1)\frac{a|Z_i|^{\beta-2}}{D} + 2\right]I_i^c\right\}$$

$$+ E\left[\left(2\beta a \frac{|Z_i|^{\beta-2}}{D}\frac{|Z_i|^\beta}{D}\right)I_i^c\right].$$

Using (A.7), one can establish that the last expression is bounded above by

$$\theta_i^2 + E[(-2(\beta-1)+2)I_i^c] + E2\beta\frac{|Z_i|^\beta}{D}I_i^c \leq \theta_i^2 + E[(4+2\beta)I_i^c]$$

(A.8)
$$\leq \theta_i^2 + 8EI_i^c.$$

We shall show under the condition in (A.2) that

(A.9) $$EI_i^c \leq C(|\theta_i|^2 + d^{\delta-1}(\log d)^{-1/2}).$$

This and (A.8) obviously establish (A.2). To prove (A.9), we shall consider two cases: (i) $0 < \beta \leq 1$ and (ii) $1 < \beta < 2$. For case (i) note that, for any $\delta > 0$, $EI_i^c$ equals

$$P(a|Z_i|^{\beta-2} \leq D) = P(D \geq a|Z_i|^{\beta-2}, |Z_i| \leq (2\ln d)^{1/2}/(1+\delta))$$
$$+ P(D \geq a|Z_i|^{\beta-2}, |Z_i| \geq (2\ln d)^{1/2}/(1+\delta)).$$

Obviously the last expression is bounded above by

(A.10) $$P(D \geq (1+\delta)^{2-\beta}dm_\beta) + P(|Z_i| \geq (2\ln d)^{1/2}/(1+\delta)).$$

Now the second term is bounded above by

(A.11) $$C(|\theta_i|^2 + (d^{1-\delta}\sqrt{\ln d})^{-1})$$

by a result in Donoho and Johnstone (1994). To find an upper bound for the first term in (A.10), note that by a simple calculus argument

$$|Z_i|^\beta \leq |\varepsilon_i|^\beta + |\theta_i|^\beta,$$

due to $0 < \beta \leq 1$. Hence the first term of (A.10) is bounded above by

$$P\left(\sum_{i=1}^d |\varepsilon_i|^\beta \geq (1+\delta)^{2-\beta} \, dm_\beta - \sum |\theta_i|^\beta\right).$$

Replacing $\sum |\theta_i|^\beta$ by the assumed upper bound in (A.2), the last displayed expression is bounded above by

(A.12) $$P\left(\sum_{i=1}^d |\varepsilon_i|^\beta \geq dm_\beta[(1+\delta)^{2-\beta} - (2-\beta)\delta^2]\right).$$

Using the inequality

$$(1+\delta)^{2-\beta} > 1 + (2-\beta)\delta,$$

one concludes that the quantity inside the bracket is bounded below by

$$1 + (2-\beta)(\delta - \delta^2) > 1.$$

Hence the probability (A.12) decays exponentially fast. This and (A.11) then establish (A.9) for $0 < \beta \leq 1$.

To complete the proof for (A.2), all we need to do is to prove (A.9) for case (ii), $1 < \beta < 2$.

Similarly to the argument for case (i), all we need to do is to show that the first term in (A.10) is bounded by (A.11). Now applying the triangle inequality,

$$D^{1/\beta} \leq \left(\sum |\varepsilon_i|^\beta\right)^{1/\beta} + \left(\sum |\theta_i|^\beta\right)^{1/\beta},$$

to the first term of (A.10) and using some straightforward algebraic manipulation, we obtain

(A.13)
$$P(D \geq (1+\delta)^{2-\beta} dm_\beta)$$
$$\leq P\left(\sum_{i=1}^d |\varepsilon_i|^\beta \geq dm_\beta \left[\left\{(1+\delta)^{(2-\beta)/\beta} - \frac{2-\beta}{2\beta}\delta^{2/\beta}\right\}^\beta\right]\right).$$

Note that

$$(1+\delta)^{(2-\beta)/\beta} \geq 1 + \frac{(2-\beta)\delta}{2\beta}$$

and consequently the quantity inside the brackets is bounded below by

$$\left[1 + \frac{2-\beta}{2\beta}(\delta - \delta^{2/\beta})\right]^\beta \geq 1 + (2-\beta)(\delta - \delta^{2/\beta})/2 > 1.$$

Now this shows that the probability on the right-hand side decreases exponentially fast. Hence inequality (A.9) is established for case (ii) and the proof for (A.2) is now completed.

To prove (A.3) for $0 \leq \beta \leq 1$, we may rewrite (A.4) as

$$E(\hat{\theta}_i - \theta_i)^2 = 1 + E(Z_i^2 - 2)I_i + E\left(|Z_i|^{2\beta-2}\left(\frac{a^2}{D^2} + \frac{2\beta a}{D^2}\right)I_i^c\right)$$
(A.14)
$$+ 2(1-\beta)E\left[\frac{|Z_i|^{\beta-2}a}{D}I_i^c\right].$$

Inequality (A.3), sharper than (A.1), can be possibly established due to the critical assumption $\beta \leq 1$, which implies that

(A.15) $$|Z_i|^{2\beta-2} \leq \left(\frac{a}{D}\right)^{-(2-2\beta)/(2-\beta)} \qquad \text{if } I_i^c = 1.$$

Note that the last term in (A.14) is obviously bounded above by $2(1-\beta)$. Furthermore, replace $|Z_i|^{2\beta-2}$ in the third term on the right-hand side of (A.14) by the upper bound in (A.15) and replace $Z_i^2$ in the second term by the upper bound

$$|Z_i|^2 < \left(\frac{a}{D}\right)^{2/(2-\beta)} \qquad \text{when } I_i = 1,$$

which follows easily for (A.5). We then obtain an upper bound for (A.14):

$$1 + E\left(\frac{a}{D}\right)^{2/(2-\beta)} + E\left[\left(\frac{a}{D}\right)^{(2\beta-2)/(2-\beta)}\left(\frac{a^2}{D^2} + 2\frac{\beta a}{D^2}\right)I_i^c\right] + 2(1-\beta)$$
$$\leq (3 - 2\beta) + CE\left(\frac{a}{D}\right)^{2/(2-\beta)}.$$

Here in the last inequality $2\beta a/D^2$ was replaced by $2\beta a^2/D^2$. To establish (A.3), obviously the only thing left to do is

(A.16) $$E\left(\frac{a}{D}\right)^{2/(2-\beta)} \leq C\ln(d).$$

This inequality can be established if we can show that

(A.17) $$E\left(\frac{d}{D}\right)^{2/(2-\beta)} \leq C,$$

since the definition of $a$ and a simple calculation show that

$$a^{2/(2-\beta)} = Ca^{2/(2-\beta)}\ln(d).$$

To prove (A.17), we apply Anderson's theorem [Anderson (1955)] which implies that $|Z_i|$ is stochastically larger than $|\varepsilon_i|$. Hence

$$E\left(\frac{d}{D}\right)^{2/(2-\beta)} \le E\Big[d/\Big(\sum |\varepsilon_i|^\beta\Big)\Big]^{2/(2-\beta)},$$

which is bounded by $A + B$. Here

$$A = E\Big[d/\Big(\sum |\varepsilon_i|^\beta\Big)\Big]^{2/(2-\beta)} I\bigg(\sum_{i=1}^d |\varepsilon_i|^\beta \le dm_\beta/2\bigg)$$

and

$$B = E\Big[d/\Big(\sum |\varepsilon_i|^\beta\Big)\Big]^{2/(2-\beta)} I\bigg(\sum_{i=1}^d |\varepsilon_i|^\beta > dm_\beta/2\bigg),$$

and as before $I(\cdot)$ denotes the indicator function.

Now $B$ is obviously bounded above by

$$(2/m_\beta)^{2/(2-\beta)} < C.$$

Also by the Cauchy–Schwarz inequality

$$A^2 \le E\Big[d/\Big(\sum |\varepsilon_i|^\beta\Big)\Big]^{4/(2-\beta)} P\bigg(\sum_{i=1}^d |\varepsilon_i|^\beta \le dm_\beta/2\bigg) < C.$$

Here the last inequality holds since the probability decays exponentially fast. This completes the proof for (A.17) and consequently for (A.3). □

Now we apply Lemma A.1 to the wavelet regression. We only prove the case $0 < \beta \le 2$. For $\beta = 0$ the proof is similar and simpler. Equivalently we shall consider the model

(A.18) $$Z_{jk} = \theta_{jk} + \varepsilon_{jk}/\sqrt{n}, \qquad k = 1, \ldots, 2^j,$$

where the $\theta_{jk}$'s are wavelet coefficients of the function $f$, and the $\varepsilon_{jk}$'s are i.i.d. standard normal random variables. For the details of reasoning supporting the above statement, see, for example, Section 9.2 of Cai (1999), following the ideas of Donoho and Johnstone (1994, 1995). Also assume that $\theta$'s live in the Besov space $B^\alpha_{p,q}(M)$ with smoothness $\alpha$ and shape parameters $p$ and $q$, that is,

(A.19) $$\sum_j 2^{jq(\alpha+1/2-1/p)}\bigg(\sum_k |\theta_{jk}|^p\bigg)^{q/p} \le M^q$$

for some positive constants $\alpha$, $p$, $q$ and $M$. The estimator $\hat{\theta}$ below for model (A.18) refers to (20) with $a$ defined in (22) and $\sigma^2 = 1/n$. For such a $\hat{\theta}$, the total risk can be decomposed into the sum of the following three quantities:

$$R_1 = \sum_{j<j_0} \sum_k E(\hat{\theta}_{jk} - \theta_{jk})^2,$$

$$R_2 = \sum_{J>j\geq j_0} \sum_k E(\hat{\theta}_{jk} - \theta_{jk})^2,$$

$$R_3 = \sum_{j\geq J} \sum_k E(\hat{\theta}_{jk} - \theta_{jk})^2,$$

where $j_0 = [\log_2(C_\delta n^{1/(2\alpha+1)})]$ and $C_\delta$ is a positive constant to be specified later. Applying (A.1) to $R_1$, which corresponds to the risk of low resolution, we establish by some simple calculation

(A.20) $$R_1 \leq C(\ln n)^{(2-\beta)/2} n^{-2\alpha/(2\alpha+1)}.$$

For $j \geq j_0$, (A.19) implies

(A.21) $$\sum_k |\theta_{jk}|^p \leq M^p 2^{-jp(\alpha+1/2-1/p)} = M^p 2^j 2^{-jp(\alpha+1/2)}.$$

Furthermore, for $p \geq \beta$

$$2^{-jp(\alpha+1/2)} \leq 2^{-j\beta(\alpha+1/2)} \leq 2^{-j_0\beta(\alpha+1/2)} = (C_\delta)^{-\beta(\alpha+1/2)} \sigma^\beta.$$

Choose $C_\delta > 0$ such that

$$M^p / C_\delta^{(1/2+\alpha)\beta} = \left(\frac{2-\beta}{2\beta}\right)^\beta \left(\frac{1}{2\alpha+1}\right)^2 m_\beta.$$

This then implies that

$$\sum_k |\theta_{jk}|^p \leq \frac{M^p}{C_\delta^{(1/2+\alpha)\beta}} 2^j \sigma^\beta$$

$$\leq \left(\frac{2-\beta}{2\beta}\right)^\beta \left(\frac{1}{2\alpha+1}\right) m_\beta 2^j \sigma^\beta,$$

satisfying the condition in (A.2) for $d = 2^j$ and $\delta = (2\alpha+1)^{-1}$.

Now for $p \geq 2$ we give an upper bound for the total risk.

From (A.2) we obtain

$$R_2 + R_3 \leq C \sum_{j \geq j_0} \sum_k \theta_{jk}^2 + o(n^{-2\alpha/(2\alpha+1)})$$

and from Hölder's inequality the first term is bounded above by

$$\sum_{j\geq j_0} 2^{j(1-2/p)} \left(\sum_k |\theta_{jk}|^p\right)^{2/p}.$$

Then inequality (A.21) gives

$$R_2 + R_3 \leq C \sum_{j \geq j_0} 2^{j(1-2/p)} 2^{-j2(\alpha+1/2-1/p)} + o(n^{-2\alpha/(2\alpha+1)})$$

$$= C \sum_{j \geq j_0} 2^{-j2\alpha} + o(n^{-2\alpha/(2\alpha+1)})$$

$$\leq C n^{-2\alpha/(2\alpha+1)}.$$

This and (A.20) imply (23) for $0 \leq \beta \leq 2$ and $p \geq 2$.

Note that for $\beta = 2$ the proof can be found in Donoho and Johnstone (1995). For $\beta \neq 2$ our proof is very different and much more involved.

To complete the proof of the theorem, we now focus on the case $0 \leq \beta \leq 2$ and $2 > p \geq \max\{1/\alpha, \beta\}$ and establish (23). We similarly decompose the risk of $\hat{\theta}$ as the sum of $R_1$, $R_2$ and $R_3$. Note that the bound for $R_1$ in (A.20) is still valid. Inequalities (A.2) and (A.3) imply

$$R_2 \leq \sum_{J \geq j \geq j_0} \sum_k \theta_{jk}^2 \wedge \frac{\log n}{n} + o\left(\frac{1}{n^{1-\delta}}\right)$$

for some constants $C > 0$. Furthermore, the inequality

$$\sum x_i \wedge A \leq A^{1-t} \sum x_i^t, \qquad x_i \geq 0, A > 0, 1 \geq t > 0,$$

implies

$$\sum_{J \geq j \geq j_0} \sum_k \theta_{jk}^2 \wedge \frac{\log n}{n} \leq \left(\frac{\log n}{n}\right)^{1-p/2} \sum_{J > j \geq j_0} \sum_k |\theta_{jk}|^p.$$

Some simple calculations, using (A.21), establish

$$R_2 \leq C \left(\frac{\log n}{n}\right)^{1-p/2} \sum_{J > j \geq j_0} 2^{-jp(\alpha+1/2-1/p)} + o(n^{-2\alpha/(2\alpha+1)})$$

(A.22)

$$\leq C (\log n)^{1-p/2} n^{-2\alpha/(2\alpha+1)}.$$

From Hölder's inequality, it can be seen that $R_3$ is bounded above by

$$\sum_{j \geq j_0} \left(\sum_k |\theta_{jk}|^p\right)^{2/p}.$$

Similarly to (A.22), we obtain the upper bound of $R_3$,

$$R_3 \leq C \sum_{j \geq J} 2^{-j2(\alpha+1/2-1/p)} = o(n^{-2\alpha/(2\alpha+1)}),$$

where $J$ is taken to be $\log_2 n$. Thus for $0 \leq \beta \leq 2$ and $2 \geq p \geq \max\{1/\alpha, \beta\}$, we have

$$\sup_{f \in B_{p,q}^\alpha} E\|\hat{\theta} - \theta\|^2 \leq C(\log n)^{1-\beta/2} n^{-2\alpha/(2\alpha+1)}.$$

**Acknowledgments.** The authors wish to thank Professor Martin Wells at Cornell University for his interesting discussions and suggestions which led to a better version of this paper. They also wish to thank Professors Lawrence D. Brown, Dipak Dey, Edward I. George, Mark G. Low and William E. Strawderman for their encouraging comments. The authors would like to thank Arletta Havlik for her enormous patience in typing this manuscript.

Department of Statistics  
Yale University  
P.O. Box 208290  
New Haven, Connecticut 06520-8290  
USA  
e-mail: [Huibin.Zhou@Yale.edu](Huibin.Zhou@Yale.edu)

Department of Mathematics  
Cornell University  
Malott Hall  
Ithaca, New York 14853  
USA  
e-mail: [hwang@math.cornell.edu](hwang@math.cornell.edu)